\documentclass[reqno,a4paper]{amsart}
\usepackage{setspace,amssymb}
\usepackage{ifpdf}
\ifpdf
 \usepackage[hyperindex,pagebackref]{hyperref}%
\else
 \expandafter\ifx\csname dvipdfm\endcsname\relax
 \usepackage[hypertex,hyperindex,pagebackref]{hyperref}
 \else
 \usepackage[dvipdfm,hyperindex,pagebackref]{hyperref}
 \fi
\fi
\allowdisplaybreaks[4]
\theoremstyle{plain}
\newtheorem{thm}{Theorem}[section]
\theoremstyle{remark}
\newtheorem{rem}{Remark}[section]
\DeclareMathOperator{\td}{d}
\newcommand{\bell}{\textup{B}}
\numberwithin{equation}{section}

\begin{document}

\title[An explicit formula for Bernoulli numbers]
{An explicit formula for Bernoulli numbers in terms of Stirling numbers of the second kind}

\author[F. Qi]{Feng Qi}
\address[Qi]{Institute of Mathematics, Henan Polytechnic University, Jiaozuo City, Henan Province, 454010, China}
\email{\href{mailto: F. Qi <qifeng618@gmail.com>}{qifeng618@gmail.com}, \href{mailto: F. Qi <qifeng618@hotmail.com>}{qifeng618@hotmail.com}, \href{mailto: F. Qi <qifeng618@qq.com>}{qifeng618@qq.com}}
\urladdr{\url{http://qifeng618.wordpress.com}}

\begin{abstract}
In the note, the author discovers an explicit formula for computing Bernoulli numbers in terms of Stirling numbers of the second kind.
\end{abstract}

\keywords{explicit formula; Bernoulli number; Stirling number of the second kind; Bell polynomial of the second kind}

\subjclass[2010]{Primary 11B68, Secondary 11B73, 11B83}

\thanks{This paper was typeset using \AmS-\LaTeX}

\maketitle

\section{Introduction}

It is well known that Bernoulli numbers $B_{k}$ for $k\ge0$ may be generated by
\begin{equation}\label{Bernumber-dfn}
\frac{x}{e^x-1}=\sum_{k=0}^\infty B_k\frac{x^k}{k!}=1-\frac{x}2+\sum_{k=1}^\infty B_{2k}\frac{x^{2k}}{(2k)!}, \quad \vert z\vert<2\pi.
\end{equation}
See~\cite[p.~48]{Comtet-Combinatorics-74}. In combinatorics, Stirling numbers $S(n,k)$ of the second kind for $n\ge k\ge1$ may be computed by
\begin{equation}\label{Stirling-Number-dfn}
S(n,k)=\frac1{k!}\sum_{\ell=1}^k(-1)^{k-\ell}\binom{k}{\ell}\ell^{n}
\end{equation}
and may be generated by
\begin{equation}\label{2stirling-gen-funct-exp}
\frac{(e^x-1)^k}{k!}=\sum_{n=k}^\infty S(n,k)\frac{x^n}{n!}, \quad k\in\mathbb{N}.
\end{equation}
See~\cite[p.~206]{Comtet-Combinatorics-74}. Bell polynomials $\bell_{n,k}(x_1,x_2,\dotsc,x_{n-k+1})$ of the second kind are defined by
\begin{equation}
\bell_{n,k}(x_1,x_2,\dotsc,x_{n-k+1})=\sum_{\substack{1\le i\le n,\ell_i\in\mathbb{N}\\ \sum_{i=1}^ni\ell_i=n\\ \sum_{i=1}^n\ell_i=k}}\frac{n!}{\prod_{i=1}^{n-k+1}\ell_i!} \prod_{i=1}^{n-k+1}\Bigl(\frac{x_i}{i!}\Bigr)^{\ell_i}
\end{equation}
for $n\ge k\ge1$, see~\cite[p.~134, Theorem~A]{Comtet-Combinatorics-74}.
\par
The aim of this note is to find an explicit formula for computing Bernoulli numbers $B_n$ in terms of Stirling numbers $S(n,k)$ of the second kind.
\par
The main results may be summarized as the following theorem.

\begin{thm}\label{Bernoulli-Stirling-thm}
For $n\ge k\ge1$, we have
\begin{align}
\bell_{n,k}(0,\overbrace{1,\dotsc,1}^{n-k})&=\sum_{i=0}^k(-1)^{i}\binom{n}{i}S(n-i,k-i), \label{B-S-0-1-value}\\
\bell_{n,k}\biggl(\frac12, \frac13,\dotsc,\frac1{n-k+2}\biggr)
&=\frac{n!}{(n+k)!}\sum_{i=0}^k(-1)^{k-i}\binom{n+k}{k-i}S(n+i,i), \label{B-S-frac-value}
\end{align}
and
\begin{equation}\label{Bernoulli-Stirling-formula}
B_n=\sum_{i=0}^n(-1)^{i}\frac{\binom{n+1}{i+1}}{\binom{n+i}{i}}S(n+i,i).
\end{equation}
\end{thm}

\section{Proof of Theorem~\ref{Bernoulli-Stirling-thm}}

In combinatorics, Fa\`a di Bruno formula may be described in terms of the second kind Bell polynomials $\bell_{n,k}(x_1,x_2,\dotsc,x_{n-k+1})$ by
\begin{equation}\label{Bruno-Bell-Polynomial}
\frac{\td^n}{\td x^n}f\circ g(x)=\sum_{k=1}^nf^{(k)}(g(x)) \bell_{n,k}\bigl(g'(x),g''(x),\dotsc,g^{(n-k+1)}(x)\bigr).
\end{equation}
See~\cite[p.~139, Theorem~C]{Comtet-Combinatorics-74}.
It is easy to see that
\begin{equation*}
\frac{x}{e^x-1}=\frac1{\int_0^1e^{xt}\td t}.
\end{equation*}
Applying in~\eqref{Bruno-Bell-Polynomial} the functions $f(y)=\frac1y$ and $y=g(x)=\int_0^1e^{xt}\td t$ results in
\begin{align*}
&\quad\frac{\td^n}{\td x^n}\biggl(\frac{x}{e^x-1}\biggr)
=\frac{\td^n}{\td x^n}\Biggl(\frac1{\int_0^1e^{xt}\td t}\Biggr)\\
&=\sum_{k=1}^n(-1)^k\frac{k!}{\bigl(\int_0^1e^{xt}\td t\bigr)^{k+1}} \bell_{n,k}\biggl(\int_0^1te^{xt}\td t, \int_0^1t^2e^{xt}\td t,\dotsc,\int_0^1t^{n-k+1}e^{xt}\td t\biggr)\\
&\to \sum_{k=1}^n(-1)^kk! \bell_{n,k}\biggl(\int_0^1t\td t, \int_0^1t^2\td t,\dotsc,\int_0^1t^{n-k+1}\td t\biggr),\quad x\to0\\
&=\sum_{k=1}^n(-1)^kk! \bell_{n,k}\biggl(\frac12, \frac13,\dotsc,\frac1{n-k+2}\biggr).
\end{align*}
On the other hand, differentiating $n$ times on both sides of~\eqref{Bernumber-dfn} leads to
\begin{equation*}
\frac{\td^n}{\td x^n}\biggl(\frac{x}{e^x-1}\biggr)=\sum_{k=n}^\infty B_k\frac{x^{k-n}}{(k-n)!}
\to B_n,\quad x\to0.
\end{equation*}
As a result, we obtain
\begin{equation}\label{Bernoulli-Bell-formula}
B_n=\sum_{k=1}^n(-1)^kk! \bell_{n,k}\biggl(\frac12, \frac13,\dotsc,\frac1{n-k+2}\biggr).
\end{equation}
\par
In~\cite[p.~113]{Comtet-Combinatorics-74}, it was listed that
\begin{equation}\label{113-final-formula}
\frac1{k!}\Biggl(\sum_{m=1}^\infty x_m\frac{t^m}{m!}\Biggr)^k =\sum_{n=k}^\infty \bell_{n,k}(x_1,x_2,\dotsc,x_{n-k+1})\frac{t^n}{n!}, \quad k\ge0.
\end{equation}
Letting $x_1=0$ and $x_m=1$ for $m\ge2$ in~\eqref{113-final-formula} and employing~\eqref{2stirling-gen-funct-exp} give
\begin{gather*}
\sum_{n=k}^\infty \bell_{n,k}(0,\overbrace{1,\dotsc,1}^{n-k})\frac{t^n}{n!}
=\frac1{k!}\Biggl(\sum_{m=2}^\infty\frac{t^m}{m!}\Biggr)^k
=\frac1{k!}(e^t-1-t)^k\\
=\frac1{k!}\sum_{i=0}^k(-1)^{k-i}\binom{k}{i}(e^t-1)^it^{k-i}
=\sum_{i=0}^k\frac{(-1)^{k-i}}{(k-i)!}\sum_{j=i}^\infty S(j,i)\frac{t^{k+j-i}}{j!}.
\end{gather*}
This implies that
\begin{multline*}
\bell_{n,k}(0,\overbrace{1,\dotsc,1}^{n-k})
=n!\sum_{i=0}^k\frac{(-1)^{k-i}}{(k-i)!}\frac{S(n-k+i,i)}{(n-k+i)!}\\
=\sum_{i=0}^k(-1)^{k-i}\binom{n}{k-i}S(n-k+i,i)
=\sum_{i=0}^k(-1)^{i}\binom{n}{i}S(n-i,k-i).
\end{multline*}
The formula~\eqref{B-S-0-1-value} follows.
\par
By virtue of
\begin{equation}
\bell_{n,k}\biggl(\frac{x_2}2, \frac{x_3}3,\dotsc,\frac{x_{n-k+2}}{n-k+2}\biggr)
=\frac{n!}{(n+k)!}\bell_{n+k,k}(0,x_2,\dotsc,x_{n+1}),
\end{equation}
see~\cite[p.~136]{Comtet-Combinatorics-74}, and the formula~\eqref{B-S-0-1-value}, we obtain
\begin{multline*}
\bell_{n,k}\biggl(\frac12, \frac13,\dotsc,\frac1{n-k+2}\biggr)
=\frac{n!}{(n+k)!}\bell_{n+k,k}(0,\overbrace{1,\dotsc,1}^n)\\
=\frac{n!}{(n+k)!}\sum_{i=0}^k(-1)^{i}\binom{n+k}{i}S(n+k-i,k-i),
\end{multline*}
from which, the formula~\eqref{B-S-frac-value} follows.
\par
Substituting~\eqref{B-S-frac-value} into~\eqref{Bernoulli-Bell-formula} leads to
\begin{gather*}
B_n=\sum_{k=1}^n\frac{k!n!}{(n+k)!}\sum_{i=0}^k(-1)^{i}\binom{n+k}{k-i}S(n+i,i)
=\sum_{k=1}^n\sum_{i=0}^k(-1)^{i}\frac{\binom{k}{i}}{\binom{n+i}{i}}S(n+i,i)\\
=\sum_{i=0}^n\frac{(-1)^{i}}{\binom{n+i}{i}}S(n+i,i)\sum_{k=i}^n\binom{k}{i}
=\sum_{i=0}^n\frac{(-1)^{i}}{\binom{n+i}{i}}\binom{n+1}{i+1}S(n+i,i),
\end{gather*}
which may be rewritten as the formula~\eqref{Bernoulli-Stirling-formula}.
The proof of Theorem~\ref{Bernoulli-Stirling-thm} is complete.

\section{Remarks}

\begin{rem}
The formula~\eqref{B-S-frac-value} may be alternatively proved as follows.
\par
Taking $x_m=\frac1{m+1}$ for all $m\in\mathbb{N}$ in~\eqref{113-final-formula} and utilizing~\eqref{2stirling-gen-funct-exp} yield
\begin{gather*}
\sum_{n=k}^\infty \bell_{n,k}\biggl(\frac12, \frac13,\dotsc,\frac1{n-k+2}\biggr)\frac{t^n}{n!}
=\frac1{k!}\Biggl[\sum_{m=1}^\infty\frac{t^m}{(m+1)!}\Biggr]^k
=\frac1{k!}\biggl(\frac{e^t-1-t}t\biggr)^k\\
=\frac1{k!}\biggl(\frac{e^t-1}t-1\biggr)^k
=\frac1{k!}\sum_{\ell=0}^k(-1)^{k-\ell}\binom{k}{\ell}\biggl(\frac{e^t-1}t\biggr)^\ell\\
=\frac1{k!}\sum_{\ell=0}^k(-1)^{k-\ell}\binom{k}{\ell} \frac{\ell!}{t^\ell} \sum_{i=\ell}^\infty S(i,\ell)\frac{t^i}{i!}
=\sum_{\ell=0}^k\frac{(-1)^{k-\ell}}{(k-\ell)!}\sum_{i=\ell}^\infty S(i,\ell)\frac{t^{i-\ell}}{i!}.
\end{gather*}
This implies that
\begin{gather*}
\bell_{n,k}\biggl(\frac12, \frac13,\dotsc,\frac1{n-k+2}\biggr)
=n!\sum_{\ell=0}^k\frac{(-1)^{k-\ell}}{(k-\ell)!(n+\ell)!}S(n+\ell,\ell).
\end{gather*}
The formula~\eqref{B-S-frac-value} follows.
\end{rem}

\begin{rem}
We collect several formulas for computing Bernoulli numbers $B_n$ as follows.
\par
In~\cite{Logan-Bell-87}, see also~\cite[pp.~559\nobreakdash--560]{GKP-Concrete-Math}, the following explicit formula for computing Bernoulli numbers $B_n$ in terms of the second kind Stirling numbers $S(n,k)$ was presented: For $n\ge1$, we have
\begin{equation}\label{Bernoulli-Stirling-eq}
B_n=\sum_{k=1}^n(-1)^k\frac{k!}{k+1}S(n,k).
\end{equation}
In~\cite[p.~1128, Corollary]{recursion}, among other things, it was found that, for $k\ge 1$,
\begin{equation}\label{Bernoulli-N-Guo-Qi-99}
B_{2k}= \frac12 - \frac1{2k+1} - 2k \sum_{i=1}^{k-1}
\frac{A_{2(k-i)}}{2(k - i) + 1},
\end{equation}
where $A_m$ is defined by
\begin{equation*}
\sum_{m=1}^nm^k=\sum_{m=0}^{k+1}A_mn^{m}.
\end{equation*}
In~\cite[Theorem~3.1]{exp-derivative-sum-Combined.tex}, it was presented that Bernoulli numbers $B_{2k}$ may be computed by
\begin{multline}\label{Bernumber-formula-eq}
B_{2k}=1+\sum_{m=1}^{2k-1}\frac{S(2k+1,m+1) S(2k,2k-m)}{\binom{2k}{m}} \\*
-\frac{2k}{2k+1}\sum_{m=1}^{2k}\frac{S(2k,m)S(2k+1,2k-m+1)}{\binom{2k}{m-1}}, \quad k\in\mathbb{N}.
\end{multline}
In~\cite[Theorem~1.4]{Tan-Cot-Bernulli-No.tex}, among other things, it was discovered that, for $k\in\mathbb{N}$,
\begin{equation}\label{Bernoulli-N-Explicit}
B_{2k}=\frac{(-1)^{k-1}k}{2^{2(k-1)}(2^{2k}-1)}\sum_{i=0}^{k-1}\sum_{\ell=0}^{k-i-1} (-1)^{i+\ell}\binom{2k}{\ell}(k-i-\ell)^{2k-1}.
\end{equation}
\end{rem}

\subsection*{Acknowledgements}
The author thanks Dr Ai-Min Xu in China and Dr Vladimir Kruchinin in Russia for their sending variants of the formula~\eqref{B-S-frac-value} through e-mail on 16 and 17 October 2013, thanks Dr Armen Bagdasaryan in Russia for his pointing out through ResearchGate on 17 October 2013, when this manuscript was completed, the formula~\eqref{B-S-frac-value} appearing in the paper~\cite{Zhang-Yang-Oxford-Taiwan-12}, and thanks Professor Juan B. Gil for his providing the formula~\eqref{B-S-frac-value} through e-mail on 24 October 2013.

\end{document}